\documentclass[12pt]{article}                
\usepackage{graphicx}
\usepackage{psfrag}
\usepackage{amsmath, amssymb}
\usepackage{mathptmx}
\usepackage{ntheorem}

\usepackage{setspace}
\makeatletter
\newsavebox\myboxA
\newsavebox\myboxB
\newlength\mylenA
\newcommand*\xoverline[2][0.75]{%
    \sbox{\myboxA}{$\m@th#2$}%
    \setbox\myboxB\null
    \ht\myboxB=\ht\myboxA%
    \dp\myboxB=\dp\myboxA%
    \wd\myboxB=#1\wd\myboxA
    \sbox\myboxB{$\m@th\overline{\copy\myboxB}$}
    \setlength\mylenA{\the\wd\myboxA}
    \addtolength\mylenA{-\the\wd\myboxB}%
    \ifdim\wd\myboxB<\wd\myboxA%
       \rlap{\hskip 0.5\mylenA\usebox\myboxB}{\usebox\myboxA}%
    \else
        \hskip -0.5\mylenA\rlap{\usebox\myboxA}{\hskip 0.5\mylenA\usebox\myboxB}%
    \fi}
\makeatother

\DeclareMathOperator{\I}{I}

\newcommand{\eps}{{\varepsilon}}
\newcommand{\ch}{{\mbox{\rm ch}}}
\newcommand{\myth}{{\mbox{\rm th}}}

\newcommand{\e}{\mathbb{E}}

\newcommand{\p}{\mathbb{P}}

\newcommand{\Natural}{\mathbb{N}}
\newcommand{\la}{\langle}
\newcommand{\ra}{\rangle}
\newcommand\qed{\hfill\hbox{\rlap{$\sqcap$}$\sqcup$}}

\newcommand{\ME}{\mbox{\it ME}}
\newcommand{\MCE}{\mbox{\it MCE}}

\newcommand{\vX}{\mathbf{X}}

\newtheorem{proposition}{Proposition}
\newtheorem{lemma}{Lemma}
\newtheorem{theorem}{Theorem}
\newtheorem{corollary}{Corollary}
\newtheorem{remark}{Remark}

\theoremstyle{nonumberplain}
\newcommand\specialref{}

\begin{document}

	\nocite{*} 
	
	\title{Disorder chaos in some diluted spin glass models}
	
	\author{Wei-Kuo Chen\thanks{School of Mathematics, University of Minnesota. Email: wkchen@umn.edu}
		\and
		Dmitry Panchenko\thanks{Department of Mathematics, University of Toronto. Email: panchenk@math.toronto.edu}}
		
\date{}
\maketitle
\begin{abstract}
We prove disorder chaos at zero temperature for three types of diluted models with large connectivity parameter: $K$-spin antiferromagnetic Ising model for even $K\geq 2$, $K$-spin spin glass model for even $K\geq 2$, and random $K$-sat model for all $K\geq 2$. We show that modifying even a small proportion of clauses results in near maximizers of the original and modified Hamiltonians being nearly orthogonal to each other with high probability. We use a standard technique of approximating diluted models by appropriate fully connected models and then apply disorder chaos results in this setting, which include both previously known results as well as new examples motivated by the random $K$-sat model.   
\end{abstract} 
\vspace{0.5cm}

\section{Introduction}

We begin by defining three types of diluted spin glass models that will be considered in this paper. Let $K\geq 2$ be a fixed integer. A random clause with $K$ variables is a random function $\theta(\sigma_{1},\ldots,\sigma_{K})$ on $\{-1,+1\}^K$. We will consider the following three examples. 

\medskip
\noindent
{\bf Example 1.} ($K$-spin antiferromagnetic model)  The function $\theta$ is given by
\begin{equation}
\theta(\sigma_1,\ldots,\sigma_K)= -\sigma_1\cdots \sigma_K,
\end{equation}
so in this case it is non-random.

\medskip
\noindent
{\bf Example 2.} ($K$-spin spin glass) The random function $\theta$ is given by
\begin{equation}
\theta(\sigma_1,\ldots,\sigma_K)= J \sigma_1\cdots \sigma_K,
\end{equation}
where $J$ is a Rademacher random variable with $\p(J=\pm 1)=1/2$.

\medskip
\noindent
{\bf Example 3.} ($K$-sat model) The random function $\theta$ is given by
\begin{equation}
\theta(\sigma_1,\ldots,\sigma_K)=-\prod_{k\leq K} \frac{1+J_k \sigma_k}{2},
\label{Ksatclause}
\end{equation}
where $(J_k)_{k\geq 1}$ are i.i.d. Rademacher random variables.

\medskip

We choose random variables $J$ to be Rademacher only for simplicity, and one can also consider other symmetric distributions. We will denote by $\theta_j$ independent copies of the function $\theta$ for various indices $j$. In the Example 1, $\theta$ is non-random so $\theta_j=\theta$, and in Examples 2 and 3,
$$
\theta_j(\sigma_1,\ldots,\sigma_K)= J_j \sigma_1\cdots \sigma_K
\,\,\mbox{ and }\,\,
\theta_j(\sigma_1,\ldots,\sigma_K)= -\prod_{k\leq K} \frac{1+J_{k,j} \sigma_k}{2}
$$ 
with i.i.d. copied $J_j$ or $J_{k,j}$ of $J$. Given a parameter $\lambda>0$, called the connectivity parameter, the Hamiltonian of the models we will be studying is defined by
\begin{equation}
H_{\lambda}(\sigma) =  \sum_{j\leq \pi(\lambda N)} \theta_j(\sigma_{i_{1,j}},\ldots, \sigma_{i_{K,j}}),
\label{Ham2}
\end{equation}
where $\pi(\lambda N)$ is a Poiss$(\lambda N)$ random variable with the mean $\lambda N$, and the coordinate indices $i_{k,j}$ are independent for different pairs $(k,j)$ and are chosen uniformly from $\{1,\ldots, N\}$. In Example 1, we will be interested only in balanced configurations belonging to the set 
\begin{equation}
V = \Bigl\{\sigma \,:\, \sum_{i\leq N} \sigma_i = 0 \mbox{ if $N\in2\Natural$ and } \sum_{i\leq N} \sigma_i = 1 \mbox{ if $N\in2\Natural+1$}\Bigr\},
\label{zeromag}
\end{equation}
while in Examples 2 and 3 the set will be the entire hypercube $V=\{-1,+1\}^N.$ Our main result below will describe a chaotic behaviour of near maximizers of this Hamiltonian under resampling a positive proportion of clauses, even arbitrarily small.

Example 1 with $K=2$ is related to the maximum bisection of the sparse Erd\"os-R\'enyi random graph $G(N,2\lambda/N)$, since we can represent a bisection by a balanced configurations $\sigma$ and write the number of edges between the two groups (up to the usual Poisson approximation) as
$$
\sum_{j\leq \pi(\lambda N)}\I(\sigma_{i_{1,j}}\not = \sigma_{i_{2,j}})=\frac{\pi(\lambda N)}{2}-\sum_{j\leq \pi(\lambda N)} \sigma_{i_{1,j}}\sigma_{i_{2,j}} = \frac{\pi(\lambda N)}{2} + H_\lambda(\sigma).
$$
Example 2 is a diluted version of the $K$-spin Sherrington-Kirkpatrick spin glass model \cite{SK}, and Example 3 corresponds to the random $K$-sat model and, in particular, finding the assignment of variables $(\sigma_i)$ that maximizes the number of satisfied clauses, which for a given clause means that at least one $\sigma_k=-J_k$ for $1\leq k\leq K$. 

For the rest of the paper we fix a correlation/proportion parameter 
\begin{equation}
t\in (0,1),
\end{equation} 
and consider two copies $H_{\lambda}^1(\sigma)$ and $H_{\lambda}^2(\sigma)$ of (\ref{Ham2}) defined in two different ways as follows.

\noindent \emph{(Resampling clauses)} In this case, the two Hamiltonians $H_{\lambda}^1(\sigma)$ and $H_{\lambda}^2(\sigma)$ will have Poiss$(t\lambda N)$ common clauses and two independent Poiss$((1-t)\lambda N)$ independently generated clauses, which means that we resample both indices of variables and random signs. For instance, in Example 1 the only randomness is in the choice of the clause indices so, for $\ell=1,2$, we define
\begin{equation}
H_{\lambda}^\ell(\sigma) = - \sum_{j\leq \pi(t\lambda N)}\sigma_{i_{1,j}}\cdots\sigma_{i_{K,j}}
-\sum_{j\leq \pi_\ell((1-t)\lambda N)}\sigma_{i_{\ell,1,j}}\cdots\sigma_{i_{\ell, K,j}},
\label{ResClause}
\end{equation}
where $\pi(t\lambda N), \pi_1((1-t)\lambda N)$ and $\pi_2((1-t)\lambda N)$ are independent Poisson random variables with the means $t\lambda N$ and $(1-t)\lambda N$ and all indices $i_{k,j}$ and $i_{\ell,k,j}$ are independent and uniform on $\{1,\ldots, N\}.$ In the Examples 2 and 3, we can also resample only random signs $J$ without resampling indices, as follows. 

\noindent \emph{(Resampling random signs)} In this case, the number of clauses, $\pi(\lambda N)$, will be the same, as well as their indices, and only random signs will be resampled. In the $K$-spin spin glass, we will replace the random variable $J_j$ in each clause by two correlated copies $J_j^1$ and $J_j^2$ such that $\e J_j^1 J_j^2 = t.$
In the $K$-sat model, we will consider two versions of the correlated copies $H_{\lambda}^1(\sigma)$ and $H_{\lambda}^2(\sigma)$.
\begin{enumerate}
\item[(a)]
In the first version, independently for each clause $\theta_j$, with probability $1-t$ we resample all random signs $J_{1,j}^1,\ldots J_{K,j}^1$ simultaneously to produce $J_{1,j}^2,\ldots J_{K,j}^2.$ In this case, $\e J_{k,j}^1 J_{k,j}^2 = t$ but the pairs $(J_{k,j}^1, J_{k,j}^2)$ within the same clause are correlated.

\item[(b)]
In the second version we will simply replace each $J_{k,j}$ in each clause by two copies such that 
$\e J_{k,j}^1 J_{k,j}^2 = t.$
In other words, we resample each $J_{k,j}^1$ with probability $1-t$ to produce $J_{k,j}^2$, but $(J_{k,j}^1, J_{k,j}^2)$ are independent for all different pairs $(k,j).$ 
\end{enumerate}
\noindent
The main result of this paper is that for all decouplings of the Hamiltonian described above and for large connectivity $\lambda$, with high probability all near maximizers are nearly orthogonal to each other. For technical reason related to Theorem \ref{Thm2} below, we will assume from now on that $K\geq 2$ is even in Examples 1 and 2, and arbitrary $K\geq 2$ in Example 3. 
\begin{theorem}\label{Thm1}
For any $\eps, t \in (0,1)$ there exists small enough $\eta>0$ such that for large enough $\lambda$ the following holds for large enough $N$ with probability at least $1-Le^{-N\eta^2/L}$: for any configurations $\sigma^1,\sigma^2\in V$ that nearly maximize the corresponding Hamiltonian,
\begin{equation}
\frac{1}{N}H_{\lambda}^\ell(\sigma^\ell) \geq \frac{1}{N}\max_{\sigma\in V} H_{\lambda}^\ell(\sigma) - \eta\sqrt{\lambda}
\,\,\mbox{ for }\,\, \ell=1,2,
\end{equation}
the overlap $R_{1,2}=N^{-1}\sum_{i\leq N} \sigma_i^1 \sigma_i^2$ between them satisfies $|R_{1,2}|\leq \eps.$
\end{theorem}
In other words, the set of near maximizers is chaotic under resampling even a small proportion of clauses. Here the constant $L$ depends only on $K$, $L=L(K)>0$, and we will see in the proof that one can take $\lambda\geq L \eta^{-6}$ also for some large enough constant $L=L(K)>0$. The definition of near maximizers is taken on the scale $\sqrt{\lambda}$, which corresponds to the scale of the maximum (see (\ref{Thm1Int0}) and (\ref{Thm1Int1}) below) except that in the random $K$-sat model one needs to subtract $-\lambda/2^K$ corresponding to the non-random part of the clauses in (\ref{Ksatclause}). 

Let us mention one standard consequence of Theorem \ref{Thm1} -- the existence of exponentially many in $N$ near maximizers of the Hamiltonian $H_\lambda$ in the Examples 1 and 2 (and in Example 3 after subtracting $-\lambda/2^K$) which are all nearly orthogonal to each other. To see this, consider the case of resampling clauses and notice that, by (\ref{Thm1Int0}), the normalized maximum over Poiss$((1-t)\lambda N)$ clauses will be of order $\sqrt{(1-t)\lambda}$ so, for $t$ close to $1$, the maximizer of the second Hamiltonian can be considered a near maximizer of the first one. Since we have exponential control of all probabilities, we can resample the Hamiltonian exponentially many times to find exponentially many near maximizers that are all nearly orthogonal to each other. 

The proof of Theorem \ref{Thm1} will be based on some known as well as new examples of disorder chaos at zero temperature for mixed $p$-spin models, and a standard comparison technique of the diluted models with large connectivity with appropriate mixed $p$-spin models using the Guerra-Toninelli interpolation from \cite{GT}. This technique has been used in various settings in \cite{DMS,Sen, PKsat, JKS} to show that the diluted $K$-spin model and $K$-spin spin glass can be approximated by the fully connected pure $K$-spin spin glass model with the Hamiltonian
\begin{equation}
H(\sigma) = \frac{1}{N^{(K-1)/2}} \sum_{1\leq i_1,\ldots,i_K\leq N}g_{i_1,\ldots,i_K} \sigma_{i_1}\cdots\sigma_{i_K}
\label{HKspin}
\end{equation}
and the $K$-sat model can be approximated by the following mixed $p$-spin Hamiltonian 
\begin{equation}
H(\sigma) = \sum_{p=1}^K  \sqrt{{K\choose p} \frac{1}{N^{p-1}}}
\sum_{1\leq i_1,\ldots,i_p\leq N}g_{i_1,\ldots,i_p} \sigma_{i_1}\cdots\sigma_{i_p},
\label{H}
\end{equation}
where the coefficients $(g_{i_1,\ldots,i_p})$ are standard Gaussian random variables independent for all $p\geq 1$ and all indices $(i_1,\ldots,i_p)$. From now on, whenever $H_\lambda$ and $H$ appear in the same formula, we assume without explicitly mentioning again that Examples 1 and 2 are paired with (\ref{HKspin}) and Example 3 is paired with (\ref{H}). For the $K$-spin antiferromagnetic Ising model from Example 1 and $K$-spin spin glass from Example 2, it was shown in \cite{DMS, Sen} that
\begin{equation}
\frac{1}{N}\e\max_{\sigma\in V} H_{\lambda}(\sigma) = \sqrt{\lambda}\frac{1}{N}\e\max_{\sigma\in V} H(\sigma) +O(\lambda^{1/3})
\label{Thm1Int0}
\end{equation}
as $\lambda\to\infty$, uniformly in $N$. For the $K$-sat model, it was shown in \cite{PKsat} that
\begin{equation}
\frac{1}{N}\e\max_{\sigma\in V} H_{\lambda}(\sigma) = -\frac{\lambda}{2^K} + \frac{\sqrt{\lambda}}{2^K}\frac{1}{N}\e\max_{\sigma\in V} H(\sigma) +O(\lambda^{1/3})
\label{Thm1Int1}
\end{equation}
as $\lambda\to\infty$, uniformly in $N$ (which was obtained earlier in \cite{LP} by the non-rigorous replica method). All the above results were stated with the error $o(\sqrt{\lambda})$ instead of $O(\lambda^{1/3})$; however, a (not so) careful inspection of the arguments as well as the proof below will make it clear that the error term can be chosen to be $O(\lambda^{1/3}).$ Upper and lower bounds on the factor in front of $\sqrt{\lambda}$ were obtained earlier in \cite{CGHS}. The main motivation for the above approximation results was due to the fact that the limit of $N^{-1}\e\max_{\sigma\in V} H(\sigma)$ on the right hand side is well known and is given by the Parisi formula \cite{Parisi79, Parisi} proved for various mixed $p$-spin models in \cite{TPF} and \cite{PUltra, PPF} and extended to zero temperature in \cite{ChenAuf}. This zero temperature formula appears below in the equation (\ref{Parisi}).

To prove Theorem \ref{Thm1}, we will connect by a similar technique two resampled diluted systems coupled by some overlap constraint to two fully connected systems correlated in an appropriate way as follows. The covariance of the Gaussian Hamiltonians above is given by
\begin{equation}
\e H(\sigma^1) H(\sigma^2) = N \xi(R_{1,2}), 
\end{equation}
where in the case of (\ref{HKspin}) and (\ref{H}) correspondingly,
\begin{equation}
\xi(s)=s^K \,\,\mbox{ and }\,\,
\xi(s) = \sum_{p=1}^K {K\choose p} s^p = (1+s)^K-1.
\label{xi}
\end{equation}
For $t\in (0,1),$ we will consider two correlated copies $H^1, H^2$ of these Hamiltonians such that either
\begin{equation}
\e H^1(\sigma^1) H^2(\sigma^2) = N t\xi(R_{1,2})
\label{Hcorr1}
\end{equation}
or
\begin{equation}
\e H^1(\sigma^1) H^2(\sigma^2) = N \xi(tR_{1,2}).
\label{Hcorr2}
\end{equation}
The first type of correlation (\ref{Hcorr1}) will be used to approximate correlated copies of the diluted models in all cases except one -- the $K$-sat model with the resampling as in (b) above, in which case the second type of correlation (\ref{Hcorr2}) will be used. Once this approximation is established, at the core of the proof of Theorem \ref{Thm1} will be the following disorder chaos result for fully connected models.
\begin{theorem}\label{Thm2} Assume that one of the following conditions is satisfied:
	\begin{itemize}
		\item[$(i)$] $H$ is given by \eqref{HKspin} with even $K$ and \eqref{Hcorr1} holds;
		\item[$(ii)$] $H$ is given by \eqref{H} for any $K\geq 2$ and either \eqref{Hcorr1} or \eqref{Hcorr2} holds.
	\end{itemize}
For any $\eps, t \in (0,1)$, there exist $\eta>0$ such that, for large enough $N$,
\begin{equation}
\label{Thm2:eq1}
\frac{1}{N}\e\max_{|R_{1,2}|>\eps}\bigl(H^1(\sigma^1)+H^2(\sigma^2)\bigr) 
\leq \frac{2}{N}\e\max_{\sigma} H(\sigma) - \eta.
\end{equation}
\end{theorem}
The case of the $K$-spin model for even $K$ was already included in the disorder chaos result in \cite{ChenHL}, which in fact covers mixed $p$-spin models with even $p$. The other two cases corresponding to the Hamiltonian (\ref{H}) will be proved in this paper and, since the proof is very similar to \cite{ChenHL}, only necessary modifications will be detailed. For a number of earlier results about disorder chaos in various settings, see \cite{Chatt08, Chatt09, Chatt14, ChenChaos0, ChenChaos, Chen15}.

\begin{remark}
	\rm The proof of Theorem \ref{Thm2} applies to general covariance structures under the following technical assumptions. Let $\xi_0(s):=t\xi(s)$ in the case (\ref{Hcorr1}) and $\xi_0(s):=\xi(ts)$ in the case (\ref{Hcorr2}). If $\xi(s)$ and $\xi_0(s)$ are both convex functions on $[-1,1]$ that satisfy $\xi_0''(s)<\xi''(|s|)$ on $[-1,1]\setminus\{0\}$ and $\zeta_+(s)$ and $\zeta_-(s)$ defined in the equation \eqref{zeta} in Lemma \ref{lem0} below are nondecreasing on $(0,1]$, then chaos in disorder \eqref{Thm2:eq1} remains true.
\end{remark}

\section{Proof of Theorem \ref{Thm1}}
We will first show how Theorem \ref{Thm2} implies Theorem \ref{Thm1}. The proofs are very similar in all cases, so we will only detail one case when only random signs are resampled and one case when clauses are resampled. 

\subsection{Resampling random signs} First, we will describe the analogue of the Guerra-Toninelli interpolation in the case of resampling of random signs in the $K$-spin spin glass and random $K$-sat model.  In the next susbsection, we will describe the interpolation in the case of resampling clauses. 

For $s\in [0,1]$, let us consider the interpolating Hamiltonian
\begin{equation}
H(s,\sigma^1,\sigma^2) = \sum_{\ell=1}^2 \bigl(\delta H_{\lambda(1-s)}^\ell(\sigma^\ell) +\sqrt{s}\beta H^\ell(\sigma^\ell)\bigr),
\label{Hampert}
\end{equation}
where the correlated Hamiltonians $H_{\lambda(1-s)}^1$ and $H_{\lambda(1-s)}^2$ are defined in the same way as $H_{\lambda}^1$ and $H_{\lambda}^2$ in the introduction, only with the connectivity parameter $\lambda$ replaced by $\lambda(1-s).$ The inverse temperature parameters $\delta>0$ and $\beta>0$ will be chosen later. Let
\begin{equation}
\varphi(s) = \frac{1}{N}\e\log \sum_{|R_{1,2}|>\eps} \exp H(s,\sigma^1,\sigma^2)
\label{varphis}
\end{equation}
be the interpolating free energy of these correlated systems coupled by the overlap constraint $|R_{1,2}|>\eps.$ We will now compute the derivative $\varphi'(s)=\mathrm{I}+\mathrm{II}$ as a sum of two terms coming from the Gaussian integration by parts and Poisson integration by parts. Let us denote by $\la\,\cdot\,\ra_s$ the average with respect to the Gibbs measure on $\{(\sigma^1,\sigma^2)\in V^2\,:\, |R_{1,2}|>\eps\}$,
\begin{equation}
G_s(\sigma^1,\sigma^2) = \frac{\exp H(s,\sigma^1,\sigma^2)}{\sum_{|R_{1,2}|>\eps} \exp H(s,\sigma^1,\sigma^2)},
\label{Gees}
\end{equation}
corresponding to the Hamiltonian $H(s,\sigma^1,\sigma^2)$. Recall that the $K$-spin spin glass and random $K$-sat model are defined on $V=\{-1,+1\}^N.$ Let us denote the i.i.d. replicas from this measure by $(\sigma^{\ell,1},\sigma^{\ell,2})$ for $\ell\geq 1$ and let us denote
$$
R_{\ell,\ell'}^{j,j'} = \frac{1}{N}\sum_{i=1}^N \sigma_i^{\ell,j}\sigma_i^{\ell',j'}.
$$
Taking the derivative in $\sqrt{s}$ in the second term in (\ref{Hampert}) and using standard Gaussian integration by parts (see e.g. \cite{SG2} or Section 1.2 in \cite{SKmodel}) gives
$$
\mathrm{I} = \frac{\beta^2}{2}\Bigl(2\xi(1)+ 2\e\bigl\la t\xi(R_{1,1}^{1,2})\bigr\ra_s
- \e\bigl\la\xi(R_{1,2}^{1,1})+\xi(R_{1,2}^{2,2})+t\xi(R_{1,2}^{1,2})+t\xi(R_{1,2}^{2,1})\bigr\ra_s\Bigr)
$$
in the case when the correlation of $H^1$ and $H^2$ is given by (\ref{Hcorr1}) and
$$
\mathrm{I} = \frac{\beta^2}{2}\Bigl(2\xi(1)+ 2\e\bigl\la \xi(tR_{1,1}^{1,2})\bigr\ra_s
- \e\bigl\la\xi(R_{1,2}^{1,1})+\xi(R_{1,2}^{2,2})+\xi(tR_{1,2}^{1,2})+\xi(tR_{1,2}^{2,1})\bigr\ra_s\Bigr)
$$
in the case when the correlation of $H^1$ and $H^2$ is given by (\ref{Hcorr2}). The rest of the calculation below is quite similar to the one in the proof of the Franz-Leone upper bound for the free energy in diluted models in \cite{FL, PT}. To differentiate with respect to $s$ in the Poisson random variable $\pi(\lambda(1-s)N)$ in the first term in (\ref{Hampert}) we use that 
$$
\frac{d}{d s}\e f(\pi(s)) = \e f(\pi(s)+1) - \e f(\pi(s))
$$
for a Poisson random variables $\pi(s)$ with the mean $s$. Therefore, the derivative of $\varphi(s)$ with respect to $s$ in $\pi(\lambda(1-s)N)$ equals
$$
\mathrm{II} =  -\lambda\Bigl(\e\log \sum_{|R_{1,2}|>\eps} \exp H^+(s,\sigma^1,\sigma^2)
- \e\log \sum_{|R_{1,2}|>\eps} \exp H(s,\sigma^1,\sigma^2)\Bigr),
$$
where $H^+(s,\sigma^1,\sigma^2)$ differs from $H(s,\sigma^1,\sigma^2)$ by the addition of one more clause in each of the correlated diluted Hamiltonians,
$$
H^+(s,\sigma^1,\sigma^2)= H^+(s,\sigma^1,\sigma^2) +\delta \theta^1(\sigma_{i_{1}}^1,\ldots, \sigma^2_{i_{K}}) +\delta \theta^2(\sigma^2_{i_{1}},\ldots, \sigma^2_{i_{K}}),
$$
and these clauses are independent of $H(s,\sigma^1,\sigma^2)$ and are given by
$$
\theta^\ell (\sigma^\ell_{i_{1}},\ldots, \sigma^\ell_{i_{K}}) = -\prod_{k\leq K} \frac{1+J_{k}^\ell \sigma^\ell_{i_k}}{2},
$$
where the random signs $J_{k}^\ell$ are correlated as in the case (a) or case (b), that is, they are resampled with probability $1-t$ either independently or simultaneously within this one clause. Clearly, we can rewrite the derivative above as
$$
\mathrm{II} =  -\lambda \e\log \bigl\la
\exp \delta \theta^1(\sigma^1_{i_{1}},\ldots, \sigma^1_{i_{K}})
\exp \delta \theta^2(\sigma^2_{i_{1}},\ldots, \sigma^2_{i_{K}})
\bigr\ra_s.
$$
Since $\theta\in\{-1,0\}$, we can write $\exp \delta \theta = 1+(1-e^{-\delta})\theta$ and
$$
\exp \delta \theta^1 \exp \delta \theta^2 = 1-(1-e^{-\delta})\Delta(\sigma^1,\sigma^2),
$$
where
$$
\Delta(\sigma^1,\sigma^2) = \prod_{k\leq K}\frac{1+J_{k}^1 \sigma^1_{i_k}}{2} + \prod_{k\leq K}\frac{1+J_{k}^2 \sigma^2_{i_k}}{2} -(1-e^{-\delta})\prod_{k\leq K} \frac{1+J_{k}^1 \sigma^1_{i_k}}{2} \cdot \frac{1+J_{k}^2 \sigma^2_{i_k}}{2}.
$$
Since $0\leq \Delta(\sigma^1,\sigma^2)\leq 1+e^{-\delta}$ and $(1-e^{-\delta})\Delta(\sigma^1,\sigma^2)\leq 1-e^{-2\delta}$, we can express the logarithm above using the Taylor series as
$$
\mathrm{II} = \lambda \sum_{n\geq 1}\frac{(1-e^{-\delta})^n}{n}
\e \bigl\la \Delta(\sigma^1,\sigma^2) \bigr\ra_s^n.
$$
Using replicas, we can represent 
$$
\e \bigl\la \Delta(\sigma^1,\sigma^2) \bigr\ra_s^n
=
\e \bigl\la \prod_{\ell\leq n}\Delta(\sigma^{\ell,1},\sigma^{\ell,2}) \bigr\ra_s
=
\e \bigl\la \e' \prod_{\ell\leq n}\Delta(\sigma^{\ell,1},\sigma^{\ell,2}) \bigr\ra_s,
$$
where $\e'$ is the expectation with respect to the randomness $J_{k}^1,J_{k}^2$ and $i_k$ of the clauses $\theta^1$ and $\theta^2$, which is independent of the randomness in $\la\,\cdot\,\ra_s$. To compute this expectation, let us first note that the expectation $\e_J$ in the random variables $J_k^\ell$ satisfies
 $$
 \e_J \prod_{k\leq K} \frac{1+J_{k}^1 \sigma^1_{i_k}}{2} \cdot \frac{1+J_{k}^1 \sigma^2_{i_k}}{2}
 =
 \e_J \prod_{k\leq K} \frac{1+J_{k}^2 \sigma^1_{i_k}}{2} \cdot \frac{1+J_{k}^2 \sigma^2_{i_k}}{2}
 =
\prod_{k\leq K} \frac{1+\sigma^1_{i_k} \sigma^2_{i_k}}{4},
$$
in the case (a) of the correlations between $J_{k}^1$ and $J_{k}^2$ we have
 $$
 \e_J \prod_{k\leq K} \frac{1+J_{k}^1 \sigma^1_{i_k}}{2} \cdot \frac{1+J_{k}^2 \sigma^2_{i_k}}{2}
=
t\prod_{k\leq K} \frac{1+\sigma^1_{i_k} \sigma^2_{i_k}}{4} + \frac{1-t}{4^K},
$$
and in the case (b) of the correlations between $J_{k}^1$ and $J_{k}^2$ we have
 $$
 \e_J \prod_{k\leq K} \frac{1+J_{k}^1 \sigma^1_{i_k}}{2} \cdot \frac{1+J_{k}^2 \sigma^2_{i_k}}{2}
=
\prod_{k\leq K} \frac{1+t\sigma^1_{i_k} \sigma^2_{i_k}}{4}.
$$
Taking expectation with respect to the random indices $i_k$, we get
 $$
 \e' \prod_{k\leq K} \frac{1+J_{k}^1 \sigma^1_{i_k}}{2} \cdot \frac{1+J_{k}^1 \sigma^2_{i_k}}{2}
 =
 \e' \prod_{k\leq K} \frac{1+J_{k}^2 \sigma^1_{i_k}}{2} \cdot \frac{1+J_{k}^2 \sigma^2_{i_k}}{2}
 =
\frac{(1+R_{1,2})^K}{4^K}
=\frac{1+\xi(R_{1,2})}{4^K},
$$
in the case (a) of the correlations between $J_{k}^1$ and $J_{k}^2$ we have
 $$
 \e' \prod_{k\leq K} \frac{1+J_{k}^1 \sigma^1_{i_k}}{2} \cdot \frac{1+J_{k}^2 \sigma^2_{i_k}}{2}
= \frac{1+t\xi(R_{1,2})}{4^K},
$$
and in the case (b) of the correlations between $J_{k}^1$ and $J_{k}^2$ we have
 $$
 \e' \prod_{k\leq K} \frac{1+J_{k}^1 \sigma^1_{i_k}}{2} \cdot \frac{1+J_{k}^2 \sigma^2_{i_k}}{2}
= \frac{1+\xi(tR_{1,2})}{4^K}.
$$
From now on we will consider the case (a) since the second case (b) is similar. Using the above formulas, let us compute $\e \la \Delta(\sigma^1,\sigma^2) \ra_s^n$ above for $n=1,2$ first. For $n=1$,
$$
\e \bigl\la \Delta(\sigma^1,\sigma^2) \bigr\ra_s
=
\frac{2}{2^K}-\frac{1-e^{-\delta}}{4^K}(1+t\e\bigl\la \xi(R_{1,1}^{1,2}) \bigr\ra_s).
$$
For $n=2$, we will separate the terms that contain the factor $(1-e^{-\delta})$ to obtain
$$
\e \bigl\la \Delta(\sigma^1,\sigma^2) \bigr\ra_s^2
=
\e\Bigl\la
\frac{1+\xi(R_{1,2}^{1,1})}{4^K}+\frac{1+\xi(R_{1,2}^{2,2})}{4^K}+\frac{1+t\xi(R_{1,2}^{1,2})}{4^K}
+\frac{1+t\xi(R_{1,2}^{2,1})}{4^K}
\Bigr\ra_s + \mathrm{III_1},
$$
where $|\mathrm{III_1}|\leq L(1-e^{-\delta}).$ Finally, using that $\Delta(\sigma^1,\sigma^2)\leq 1+e^{-\delta}$, we can bound the sum over $n\geq 3$ in $\mathrm{II}$ above by
$$
\mathrm{III}_2 :=  \lambda \sum_{n\geq 3}\frac{(1-e^{-\delta})^n}{n}
\e \bigl\la \Delta(\sigma^1,\sigma^2) \bigr\ra_s^n
\leq
\lambda \sum_{n\geq 3}\frac{(1-e^{-2\delta})^n}{n}
\leq L \lambda(1-e^{-2\delta})^3.
$$
Plugging all these expressions back into $\mathrm{II}$,
\begin{align*}
\mathrm{II} = &\,\,
\frac{2\lambda(1-e^{-\delta})}{2^K}-\frac{\lambda(1-e^{-\delta})^2}{4^K}\bigl(1+t\e\bigl\la \xi(R_{1,1}^{1,2})\bigr\ra_s\bigr)
\\
&+\frac{4\lambda(1-e^{-\delta})^2}{2\cdot 4^K}
+\frac{\lambda(1-e^{-\delta})^2}{2\cdot 4^K} \e\bigl\la \xi(R_{1,2}^{1,1})+\xi(R_{1,2}^{2,2})+t\xi(R_{1,2}^{1,2})+t\xi(R_{1,2}^{2,1})\bigr\ra_s + \mathrm{III},
\end{align*}
where $|\mathrm{III}|\leq L\lambda(1-e^{-2\delta})^3.$ 

Next, given $\lambda$ and $\delta$, we are going to make the following choice of $\beta$,
\begin{equation}
\frac{\beta^2}{2} = \frac{{\lambda}(1-e^{-\delta})^2}{2\cdot 4^K}, \,\mbox{ or }\,
\beta = \frac{\sqrt{\lambda}(1-e^{-\delta})}{2^K}.
\end{equation}
With this choice, all the terms containing overlaps $R_{\ell,\ell'}^{j,j'}$ cancel out and, since $2\xi(1)=2\cdot 2^K-2$, 
$$
\varphi'(s) = \mathrm{I} + \mathrm{II} = \frac{2\lambda(1-e^{-\delta})}{2^K}+\frac{{2\lambda}(1-e^{-\delta})^2}{2\cdot 2^K}+\mathrm{III},
$$
where $|\mathrm{III}|\leq L\lambda(1-e^{-2\delta})^3\leq L\lambda \delta^3.$ Recall that we are interested in the regime when $\lambda$ is large enough, and below we will take $\delta=\lambda^{-1/3}$. By Taylor's expansion, 
$$
\frac{2\lambda(1-e^{-\delta})}{2^K}+\frac{{2\lambda}(1-e^{-\delta})^2}{2\cdot 2^K} 
= 
\frac{2\lambda\delta}{2^K} + O(\lambda \delta^3),
$$
so $\varphi'(s) =2\lambda\delta/2^K + O(\lambda \delta^3).$ Integrating between $0$ and $1$, we get
$$
\Bigl| \varphi(0) + \frac{2\lambda\delta}{2^K} - \varphi(1)\Bigr| = O(\lambda \delta^3)
$$
and, dividing both sides by $\delta$,
$$
\Bigl| \frac{1}{\delta} \varphi(0) + \frac{2\lambda}{2^K} - \frac{1}{\delta}\varphi(1)\Bigr| = O(\lambda \delta^2).
$$
Using elementary estimates
$$
\frac{1}{N}\e\max_{|R_{1,2}|>\eps} H(s,\sigma^1,\sigma^2) \leq \varphi(s) \leq 2\log 2 + \frac{1}{N}\e\max_{|R_{1,2}|>\eps} H(s,\sigma^1,\sigma^2)
$$
for $s=0$ and $s=1$, we get
$$
\Bigl|\frac{1}{N}\e\max_{|R_{1,2}|>\eps}\bigl(H_\lambda^1(\sigma^1)+H_\lambda^2(\sigma^2)\bigr) -  \frac{1}{\delta} \varphi(0)\Bigr| \leq \frac{2\log 2}{\delta}
$$
and
$$
\Bigl| \frac{\beta}{\delta} \frac{1}{N}\e\max_{|R_{1,2}|>\eps}\bigl(H^1(\sigma^1)+H^2(\sigma^2)\bigr) 
- \frac{1}{\delta} \varphi(1) \Bigr| \leq \frac{2\log 2}{\delta}.
$$
By Taylor's expansion and our choice of $\beta$ above,
$$
\frac{\beta}{\delta}= \frac{\sqrt{\lambda}(1-e^{-\delta})}{2^K\delta} = \frac{\sqrt{\lambda}}{2^K} +O(\sqrt{\lambda\delta^2})
$$
and, therefore,
\begin{align*}
\frac{1}{N}\e\max_{|R_{1,2}|>\eps}\bigl(H_\lambda^1(\sigma^1)+H_\lambda^2(\sigma^2)\bigr)
= &\,\,
- \frac{2\lambda}{2^K}
+ \frac{\sqrt{\lambda}}{2^K}  \frac{1}{N}\e\max_{|R_{1,2}|>\eps}\bigl(H^1(\sigma^1)+H^2(\sigma^2)\bigr) 
\\
& + O\Bigl(\frac{1}{\delta}+ \lambda \delta^2+\sqrt{\lambda\delta^2}\Bigr).
\end{align*}
With the choice of $\delta = \lambda^{-1/3}$, the error term here is $O(\lambda^{1/3}).$ By Theorem \ref{Thm2}, we conclude that there exists $\eta>0$ such that
$$
\frac{1}{N}\e\max_{|R_{1,2}|>\eps}\bigl(H_\lambda^1(\sigma^1)+H_\lambda^2(\sigma^2)\bigr)
\leq
- \frac{2\lambda}{2^K}
+ \frac{2\sqrt{\lambda}}{2^K}  \frac{1}{N}\e\max_{\sigma} H(\sigma) -  \frac{2\sqrt{\lambda}\eta}{2^K} 
+ L \lambda^{1/3}
$$
and, using (\ref{Thm1Int1}), we get
$$
\frac{1}{N}\e\max_{|R_{1,2}|>\eps}\bigl(H_\lambda^1(\sigma^1)+H_\lambda^2(\sigma^2)\bigr)
\leq
\frac{2}{N}\e\max_{\sigma}H_\lambda(\sigma) -  \frac{2\sqrt{\lambda}\eta}{2^K} 
+ L \lambda^{1/3}.
$$
For $\lambda\geq L\eta^{-6}$ for large enough constant $L=L(K)$, this implies that
$$
\frac{1}{N}\e\max_{|R_{1,2}|>\eps}\bigl(H_\lambda^1(\sigma^1)+H_\lambda^2(\sigma^2)\bigr)
\leq
\frac{2}{N}\e\max_{\sigma}H_\lambda(\sigma) -  \frac{\sqrt{\lambda}\eta}{L}.
$$
By Azuma's inequality, this implies that (increasing value of  the constant $L=L(K)$)
$$
\frac{1}{N}\max_{|R_{1,2}|>\eps}\bigl(H_\lambda^1(\sigma^1)+H_\lambda^2(\sigma^2)\bigr)
\leq
\frac{1}{N}\bigl(\max_{\sigma}H_\lambda^1(\sigma)+\max_{\sigma}H_\lambda^2(\sigma)\bigr) -  \frac{\sqrt{\lambda}\eta}{L}
$$
with probability at least $1-Le^{-N\eta^2/L}.$ On this event, the existence of $\sigma^1,\sigma^2$ such that
$$
\frac{1}{N}H_{\lambda}^\ell(\sigma^\ell) \geq \frac{1}{N}\max_{\sigma} H_{\lambda}^\ell(\sigma) - \frac{\sqrt{\lambda}\eta}{3L}
$$
and such that $|R_{1,2}|>\eps$ would, obviously, lead to contradiction, so the proof is finished in the case of resampling random signs.

\subsection{Resampling clauses} Next we describe the interpolation in the case of resampling clauses. We only consider the first example of antiferromagnetic $K$-spin model. Under resampling of clauses, the Hamiltonians are defined in (\ref{ResClause}). Now, for $s\in [0,1]$, we will replace these Hamiltonians by
$$
H_{\lambda,s}^\ell(\sigma) = - \sum_{j\leq \pi((1-s)t\lambda N)}\sigma_{i_{1,j}}\cdots\sigma_{i_{K,j}}
-\sum_{j\leq \pi_\ell((1-s)(1-t)\lambda N)}\sigma_{i_{\ell,1,j}}\cdots\sigma_{i_{\ell, K,j}}
$$
and consider the interpolating Hamiltonian
\begin{equation}
H(s,\sigma^1,\sigma^2) = \sum_{\ell=1}^2 \bigl(\delta H_{\lambda,s}^\ell(\sigma^\ell) +\sqrt{s}\beta H^\ell(\sigma^\ell)\bigr).
\label{HampertC}
\end{equation}
Here, $H^1$ and $H^2$ are pure fully connected $K$-spin Hamiltonians (\ref{HKspin}) with the correlation in (\ref{Hcorr1}). We will use the same notation $\varphi(s)$ and $G_s$ as in (\ref{varphis}) and (\ref{Gees}). Recall that in this model, the pairs of configurations as well as the Gibbs measure are defined on $\{(\sigma^1,\sigma^2)\in V^2\,:\, |R_{1,2}|>\eps\}$, where $V$ is the set of configurations in (\ref{zeromag}) with zero magnetization. 

As in the previous section, the derivative of the interpolating free energy $\varphi(s)$ in (\ref{varphis}) with respect to $\sqrt{s}$ in front of the second term gives
$$
\mathrm{I} = \frac{\beta^2}{2}\Bigl(2+ 2\e\bigl\la t(R_{1,1}^{1,2})^K\bigr\ra_s
- \e\bigl\la(R_{1,2}^{1,1})^K+(R_{1,2}^{2,2})^K+t(R_{1,2}^{1,2})^K+t(R_{1,2}^{2,1})^K\bigr\ra_s\Bigr).
$$
On the other hand, the derivative with respect to $s$ in the Poisson random variables will now be applied to three different terms $\pi((1-s)t\lambda N), \pi_1((1-s)(1-t)\lambda N)$ and $\pi_2((1-s)(1-t)\lambda N)$ and, similarly to the computation above,
\begin{align*}
\mathrm{II} = &\,\,  -t\lambda \e\log \bigl\la
\exp\bigl(- \delta \sigma^1_{i_{1}}\cdots\sigma^1_{i_{K}} - \delta \sigma^2_{i_{1}}\cdots \sigma^2_{i_{K}} \bigr)
\bigr\ra_s
\\
&
-(1-t)\lambda \e\log \bigl\la
\exp\bigl(- \delta \sigma^1_{i_{1}}\cdots\sigma^1_{i_{K}}\bigr)
\bigr\ra_s
-(1-t)\lambda \e\log \bigl\la
\exp\bigl( - \delta \sigma^2_{i_{1}}\cdots \sigma^2_{i_{K}} \bigr)
\bigr\ra_s.
\end{align*}
Since the product of spins takes values $\pm 1$, we can represent
$$
\exp\bigl( - \delta \sigma^\ell_{i_{1}}\cdots \sigma^\ell_{i_{K}} \bigr)
=
\ch(\delta)\bigl(1-\myth(\delta)\sigma^\ell_{i_{1}}\cdots \sigma^\ell_{i_{K}}\bigr)
$$
and rewrite their product as
$$
\exp\bigl(- \delta \sigma^1_{i_{1}}\cdots\sigma^1_{i_{K}} - \delta \sigma^2_{i_{1}}\cdots \sigma^2_{i_{K}} \bigr)
= 
\ch(\delta)^2\bigl(1-\myth(\delta)\Delta(\sigma^1,\sigma^2)\bigr)
$$
with the notation
$$
\Delta(\sigma^1,\sigma^2) = \sigma^1_{i_{1}}\cdots\sigma^1_{i_{K}} + \sigma^2_{i_{1}}\cdots \sigma^2_{i_{K}}
-\myth(\delta)\sigma^1_{i_{1}}\cdots\sigma^1_{i_{K}}\sigma^2_{i_{1}}\cdots \sigma^2_{i_{K}}.
$$
Expressing the logarithm by its Taylor series, we can rewrite
$$
\mathrm{II}=
-2\lambda\log\ch\delta
+\lambda \sum_{n\geq 1}\frac{\myth(\delta)^n}{n}
\Bigl(
t\e \bigl\la \Delta(\sigma^1,\sigma^2) \bigr\ra_s^n
+(1-t)\sum_{\ell=1}^2\e \bigl\la \sigma^\ell_{i_{1}}\cdots\sigma^\ell_{i_{K}} \bigr\ra_s^n
\Bigr).
$$
Let us recall that, in this example, we restrict configurations to the set $V$ with zero magnetization, 
$m(\sigma^\ell)=N^{-1}\sum_{i\leq N} \sigma_i^\ell =0$. For odd $N$, magnetization equals $1/N$ which for simplicity of notation we denote by $0$. Averaging in the random signs, we get
$$
\e \bigl\la \sigma^\ell_{i_{1}}\cdots\sigma^\ell_{i_{K}} \bigr\ra_s = 
\e \bigl\la m(\sigma^\ell)^K \bigr\ra_s = 0
$$
and
$$
\e \bigl\la \Delta(\sigma^1,\sigma^2) \bigr\ra_s
=
-\myth(\delta) \e \bigl\la (R_{1,1}^{1,2})^K \bigr\ra_s.
$$
This is the only place where we used the assumption that magnetization is zero. Using replicas as above, averaging in random signs and collecting all the terms of the order $O(\lambda\myth(\delta)^3)$ into error term, it is easy to check that 
\begin{align*}
\mathrm{II}=&\,\,
-2\lambda\log\ch\delta
-t\lambda \myth(\delta)^2 \e \bigl\la (R_{1,1}^{1,2})^K \bigr\ra_s
\\ &
+\frac{\lambda \myth(\delta)^2}{2}\e\bigl\la(R_{1,2}^{1,1})^K+(R_{1,2}^{2,2})^K+t(R_{1,2}^{1,2})^K+t(R_{1,2}^{2,1})^K\bigr\ra_s + O(\lambda\delta^3).
\end{align*}
If we now take $\beta = \sqrt{\lambda}\myth(\delta)$ then all the terms containing overlaps cancel out and we get
$$
\varphi'(s)= \mathrm{I}+\mathrm{II}= -2\lambda\log\ch\delta + \lambda\myth(\delta)^2+O(\lambda\delta^3).
$$
One can check that $-2\log\ch\delta + \myth(\delta)^2=O(\delta^4)$ as $\delta\to 0$ and, therefore, $\varphi'(s)=O(\lambda\delta^3).$ By Taylor's expansion and our choice of $\beta$,
$$
\frac{\beta}{\delta}= \frac{\sqrt{\lambda}\myth(\delta)}{\delta} = \sqrt{\lambda} +O(\sqrt{\lambda}\delta^2)
$$
and, using elementary estimates connecting the free energy and maximum as above, one can check that the obtained control of the derivative implies that
$$
\frac{1}{N}\e\max_{|R_{1,2}|>\eps}\bigl(H_\lambda^1(\sigma^1)+H_\lambda^2(\sigma^2)\bigr)
= 
\sqrt{\lambda} \frac{1}{N}\e\max_{|R_{1,2}|>\eps}\bigl(H^1(\sigma^1)+H^2(\sigma^2)\bigr) 
+O(\lambda^{1/3}),
$$
if we again take $\delta=\lambda^{-1/3}$. The maximum on the right hand side is taken over $\sigma^1,\sigma^2 \in V$ with zero magnetization and we can bound it from above by removing this magnetization constraint and then applying Theorem \ref{Thm2} to the right hand side,
$$
\frac{1}{N}\e\max_{|R_{1,2}|>\eps}\bigl(H^1(\sigma^1)+H^2(\sigma^2)\bigr) 
\leq \frac{2}{N}\e\max_{\sigma} H(\sigma) - \eta,
$$
where the maximum on the right hand side is now taken over all $\sigma\in\{-1,+1\}^N$. However, since the ground state energy over the whole space is essentially the same as over subset $V$ of configurations with zero magnetization as $N\to\infty$ (see Lemma \ref{lem1} below), by reducing $\eta$ we can replace the maximum above by the one over $\sigma\in V.$ Then, using (\ref{Thm1Int0}), the proof is finished in exactly the same way as above.
\qed

\begin{lemma}
	\label{lem1}
	We have that
	\begin{align*}
	\lim_{N\rightarrow\infty}\e\max_{\sigma}\frac{H(\sigma)}{N}=\lim_{N\rightarrow\infty}\e\max_{\sigma\in V}\frac{H(\sigma)}{N}
	\end{align*}
\end{lemma}
\textbf{Proof.} Obviously the left-hand side is no less than the right-hand side. It remains to show the reverse inequality. First, we note that it is already known from the proof of \cite[Proposition 9]{ChenHL} that for any $\varepsilon>0,$
\begin{align}\label{lem2:proof:eq0}
\lim_{N\rightarrow\infty}\e\max_{|m(\sigma)|<\varepsilon}\frac{H(\sigma)}{N}=\lim_{N\rightarrow\infty}\e\max_{\sigma}\frac{H(\sigma)}{N},
\end{align}
where $m(\sigma):=N^{-1}\sum_{i=1}^N\sigma_i$ is the magnetization of $\sigma.$
Let $\varepsilon>0$ be fixed. For any $\sigma$ satisfying $|m(\sigma)|<\varepsilon$, we can find $\pi(\sigma)\in \{-1,1\}^N$ with $m(\pi(\sigma))\in V$ such that the Hamming distance between $\sigma$ and $\pi(\sigma)$ satisfies
\begin{align}\label{lem2:proof:eq1}
d\bigl(\sigma,\pi(\sigma)\bigr)&=\frac{1}{N}\sum_{i=1}^NI\bigl(\sigma_i\neq \pi(\sigma)_i\bigr)< \varepsilon
\end{align} 
provided that $N\varepsilon>1.$ Indeed, let $A=\{i:\sigma_i=1\}$ and $B=\{i:\sigma_i=-1\}$ and assume $|A|<|B|$ ($|A|\geq |B|$ is similar). If $N$ is even, we consider a partition $\{B_1,B_2,B_3\}$ of $B$ with $|B_1|=|B_2|=(|B|-|A|)/2$. Set $\pi(\sigma)$ as
\begin{align*}
\pi(\sigma)_i&=\left\{
\begin{array}{ll}
1,&\mbox{if $i\in A\cup B_1$},\\
-1,&\mbox{if $i\in B_2\cup B_3$}.
\end{array}
\right.
\end{align*}
Then clearly $m(\pi(\sigma))=0$ and 
\begin{align*}
d\bigl(\sigma,\pi(\sigma)\bigr)&=\frac{1}{N}|B_1|=\frac{|B|-|A|}{2N}=\frac{|m(\sigma)|}{2}<\frac{\varepsilon}{2}.
\end{align*}
Similarly, if $N$ is odd, then we consider a partition $\{B_1,B_2,B_3,B_4\}$ of $B$ with $|B_1|=|B_2|=(|B|-|A|-1)/2$ and $|B_4|=1$. Set $\pi(\sigma)$ by
\begin{align*}
\pi(\sigma)_i&=\left\{
\begin{array}{ll}
1,&\mbox{if $i\in A\cup B_1\cup B_4$},\\
-1,&\mbox{if $i\in B_2\cup B_3$}.
\end{array}
\right.
\end{align*}
Then $m(\pi(\sigma))=1/N$ and if $N\varepsilon>1,$
\begin{align*}
d\bigl(\sigma,\pi(\sigma)\bigr)&=\frac{1}{N}(|B_1|+|B_4|)=\frac{|B|-|A|+1}{2N}<\frac{\varepsilon}{2}+\frac{1}{2N}.
\end{align*}
All these imply that $\pi(\sigma)\in V$ and \eqref{lem2:proof:eq1} holds. Now, for any $N\varepsilon>1$, if $|m(\sigma)|<\varepsilon$, then 
\begin{align*}
\e\bigl(H(\sigma)-H(\pi(\sigma))\bigr)^2&=2N\bigl(\xi(1)-\xi\bigl(R(\sigma,\pi(\sigma))\bigr)\\
&\leq 2N\xi'(1)d(\sigma,\pi(\sigma))\\
&<2N\xi'(1)\varepsilon.
\end{align*}
We apply the Slepian inequality (see e.g. \cite{LT}) to get
\begin{align*}
\e\max_{|m(\sigma)|<\varepsilon}|H(\sigma)-H(\pi(\sigma))|\leq N\sqrt{2\varepsilon\xi'(1)\log 2}.
\end{align*}
Consequently, 
\begin{align*}
\e \max_{|m(\sigma)|<\varepsilon}\frac{H(\sigma)}{N}\leq \e \max_{\sigma\in V}\frac{H(\sigma)}{N}+\sqrt{2\varepsilon\xi'(1)\log 2}.
\end{align*}
From \eqref{lem2:proof:eq0}, sending $N\rightarrow \infty$ and then $\varepsilon \downarrow 0$ gives
\begin{align*}
\lim_{N\rightarrow\infty}\e\max_{\sigma}\frac{H(\sigma)}{N}\leq \lim_{N\rightarrow\infty}\e\max_{\sigma\in V}\frac{H(\sigma)}{N}
\end{align*} 
and this completes our proof.
\qed

\section{Proof of Theorem \ref{Thm2}}

In section, we establish Theorem \ref{Thm2} assuming the case $(ii)$. In this case $H$ is associated to
\begin{align}\label{eq-1}
\xi(s)=(1+s)^{K}-1
\end{align}
for $K\geq 2.$
For notational convenience, we denote the covariance between $H^1$ and $H^2$ by $\xi_0$, 
$$
\e H^1(\sigma)H^2(\sigma^2)=\xi_0(R_{1,2}).
$$
Throughout the remainder of this section, we assume that $t\in(0,1)$ and $H^1,H^2$ have the covariance structure \eqref{Hcorr1} or \eqref{Hcorr2}, that is, either
\begin{align}\label{eq-2}
\xi_0(s)&=t\xi(s)
\end{align}
or 
\begin{align}
\label{eq-3}
\xi_0(s)&=\xi(ts).
\end{align}
Recently, chaos in disorder for the ground energy was proved in \cite{ChenHL} for the case when $\xi$ is an infinite series of even mixtures and \eqref{eq-2} holds. The major difference between \cite{ChenHL} and our current situation is that $\xi$ now includes odd $p$-spin interactions and $H^1$ and $H^2$ possess a new type of covariance structure \eqref{eq-3}. Our proof of Theorem \ref{Thm2} is essentially based on \cite{C15}, where the author established chaos in disorder at positive temperature using the Guerra replica symmetry breaking (RSB) bound \cite{Guerra} and its two-dimensional extension (in the spirit of Talagrand \cite{TPF, TalUltra}) under general assumptions on $\xi$ and $\xi_0$. While several properties of the Parisi measure $\gamma_P(ds)$ and the Parisi PDE $\Phi_\gamma$ established in \cite{ChenHL} will be used, our proof will follow closely the one in \cite{C15} only now at zero temperature.

\subsection{Parisi formula and the Guerra-Talagrand RSB bound}
The first key ingredient is played by the Parisi formula for the maximum Hamiltonian in \cite{ChenAuf}. Let $\mathcal{U}$ be the set of all nonnegative and nondecreasing right-continuous functions $\gamma$ on $[0,1]$ such that 
$$\int_0^1\gamma(s)ds<\infty.$$ We equip the space $\mathcal{U}$ with $L^1(dx)$ norm. Define the Parisi functional on $\mathcal{U}$ by
\begin{align}\label{pf}
\mathcal{P}(\gamma)=\Phi_\gamma(0,0)-\frac{1}{2}\int_0^1s\xi''(s)\gamma(s)ds,
\end{align}
where $\Phi_\gamma(0,0)$ is defined through the weak solution of the Parisi PDE with boundary condition $\Phi_\gamma(1,x)=|x|$,
\begin{align}
\label{pde}
\partial_s\Phi_\gamma(s,x)&=-\frac{\xi''(s)}{2}\bigl(\partial_{xx}\Phi_{\gamma}(s,x)+\gamma(s)\bigl(\partial_x\Phi_\gamma(s,x)\bigr)^2\bigr)
\end{align}
for $(s,x)\in[0,1)\times\mathbb{R}.$ The existence and regularity properties of $\Phi_\gamma$ can be found in \cite{ChenHL}. The Parisi formula for the maximum energy states that, for $\ME_N=N^{-1}\max_\sigma H(\sigma),$
\begin{align}
\label{Parisi}
\ME&:=\lim_{N\rightarrow\infty}\ME_N=\inf_{\gamma\in \mathcal{U}}\mathcal{P}(\gamma)\,\,a.s.
\end{align}
Here the minimizer on the right-hand side exists and is unique, see \cite{ChenHL}. Denote this minimizer by $\gamma_P$ and call $\gamma_P(ds)$ the Parisi measure. We mention that \eqref{Parisi} was indeed established for general mixtures, see \cite{ChenAuf}.

Next, we state Guerra-Talagrand's RSB bound. For any measurable $A\in [-1,1]$, consider the normalized maximum
\begin{align*}
\MCE_N(A)&=\frac{1}{N}\max_{(\sigma^1,\sigma^2)\in A}\bigl(H^1(\sigma^1)+H^2(\sigma^2)\bigr).
\end{align*}
If $A$ contains only one point $q$, we simply write $\MCE_N(A)=\MCE_N(q).$ Let $S_N=\{R_{1,2}:\sigma^1,\sigma^2\in \{-1,+1\}^N\}$ be the set of possible overlap value for a given $N.$ Let $q\in [-1,1]$ be fixed. Denote by $\iota=1$ if $q\geq 0$ and $\iota=-1$ if $q<0.$ Define a matrix-valued function $T$ on $[0,1]$ by
\begin{align}\label{Tees}
T(s)&=\left\{
\begin{array}{ll}
\left[
\begin{array}{cc}
\xi''(s)&\iota\xi_0''(\iota s)\\
\iota \xi_0''(\iota s)&\xi''(s)
\end{array}
\right],&\mbox{if $s\in[0,|q|)$},\\
\\
\left[
\begin{array}{cc}
\xi''(s)&0\\
0&\xi''(s)
\end{array}
\right],&\mbox{if $s\in[|q|,1]$}.
\end{array}
\right.
\end{align}
For any $\gamma \in \mathcal{U}$, consider the weak solution $\Psi_\gamma(s,\mathbf{x})$ on $[0,1]\times\mathbb{R}^2$ of
\begin{align*}
\partial_s\Psi_{\gamma}(s,\mathbf{x})&=-\frac{1}{2}\bigl(\la T(s),\triangledown^2 \Psi_\gamma(s,\mathbf{x})\ra+\gamma(s)\la T(s)\triangledown \Psi_\gamma(s,\mathbf{x}),\triangledown \Psi_\gamma(s,\mathbf{x})\ra\bigr)
\end{align*}
with boundary condition $\Psi_\gamma(1,\mathbf{x})=|x_1|+|x_2|.$ Here the existence of this PDE and its regularity can be argued in a similar way as those for $\Phi_\gamma$ appearing in \cite[Appendix]{ChenHL}. Define
\begin{align}\label{eq4}
\mathcal{T}_q(\gamma)=\Psi_{\gamma}(0,0,0)-\frac{1}{2}\Bigl(\int_0^1s\xi''(s)\gamma(s)ds+\int_0^{|q|}s\xi_0''(\iota s)\gamma(s)ds\Bigr)
\end{align}
for all $\gamma\in \mathcal{U}.$ For any $q\in S_N,$ the Guerra-Talagrand RSB bound for the expected value of the normalized maximum defined above is given by 
\begin{align}\label{GT}
\e \MCE_N(q)&\leq \mathcal{T}_q(\gamma). 
\end{align}
This inequality is obtained from the usual Guerra-Talagrand upper bound at positive temperature, which holds because of the convexity of $\xi$ and $\xi_0$ on $[-1,1]$, by taking zero temperature limit with the same rescaling of the functional order parameter as in one dimensional case explained in Lemma 2 in \cite{ChenAuf}. In other words, this is a standard two-dimensional analogue of Lemma 2 in \cite{ChenAuf}.

From \eqref{Parisi} and \eqref{GT}, the proof of Theorem \ref{Thm2} relies on finding $\gamma\in \mathcal{U}$ such that
$$
\mathcal{T}_q(\gamma)<2\mathcal{P}(\gamma_P),
$$ 
whenever $|q|>\eps.$ From \eqref{pf} and \eqref{eq4}, one way of attaining this is to find $\gamma\in \mathcal{U}$ such that the following two conditions are satisfied:
\begin{itemize}
	\item[$(i)$] $\int_0^1s\xi''(s)\gamma(s)ds+\int_0^{|q|}s\xi_0''(\iota s)\gamma(s)ds=\int_0^1s\xi''(s)\gamma_P(s)ds;$
	\item[$(ii)$] $\Psi_{\gamma}(0,0,0)<2\Phi_{\gamma_P}(0,0).$
\end{itemize}
To obtain $(i)$, we choose 
\begin{align}\label{eq5}
\gamma_q(s)&:=\left\{
\begin{array}{ll}
\frac{\xi''(s)\gamma_P(s)}{\xi''(s)+\xi_0''(\iota s)},&\mbox{if $s\in[0,|q|)$},\\
\\
\gamma_P(s),&\mbox{if $s\in[|q|,1]$},
\end{array}\right.
\end{align}
and an algebraic manipulation gives
\begin{align}\label{eq7}
\int_0^1s\xi''(s)\gamma_q(s)ds+\int_0^{|q|}s\xi_0''(\iota s)\gamma_q(s)ds&=\int_0^1s\xi''(s)\gamma_P(s)ds.
\end{align}
One technical condition we need here is the requirement that $\gamma_q$ must lie in $\mathcal{U}.$ The lemma below justifies this condition.

\begin{lemma}\label{lem0} Consider $\xi$ defined through \eqref{eq-1} and $\xi_0$ defined by either \eqref{eq-2} or \eqref{eq-3}.
	We have that $\xi_0''(s)<\xi''(|s|)$ for all $s\in[-1,1]\setminus\{0\}$ and $\zeta_+$ and $\zeta_-$ are nondecreasing on $(0,1]$, where for $s\in(0,1]$,
	\begin{align}
	\begin{split}
	\label{zeta}
	\zeta_+(s)&:=\frac{\xi''(s)}{\xi''(s)+\xi_0''(s)},\\
	\zeta_-(s)&:=\frac{\xi''(s)}{\xi''(s)+\xi_0''(-s)}.
	\end{split}
	\end{align}
\end{lemma} 
Since $\zeta_+$ and $\zeta_-$ are nondecreasing and $\xi''$ and $\xi_0''$ are nonnegative, we see that $\gamma_q$ is nonnegative and nondecreasing with $\lim_{s\rightarrow |q|-}\gamma_q(s)\leq \gamma_P(|q|).$ Thus, $\gamma_q\in \mathcal{U}.$ The fact that $\xi_0''(s)<\xi''(|s|)$ is not needed for this statement, but will be used in the subsequent sections.

	{\noindent \bf Proof of Lemma \ref{lem0}.} Compute directly
	\begin{align*}
	\zeta_+'(s) &=\frac{\xi'''(s) \xi_0''(s) -\xi''(s) \xi_0'''(s) }{(\xi''(s) +\xi_0''(s) )^2}
	\end{align*}
	and
	\begin{align*}
	\zeta_-'(s) &=\frac{\xi'''(s) \xi_0''(-s) +\xi''(s) \xi_0'''(-s) }{(\xi''(s) +\xi_0''(s) )^2}.
	\end{align*}
	If \eqref{eq-2} holds, then
	\begin{align*}
	\xi'''(s) \xi_0''(s) -\xi''(s) \xi_0'''(s) =tK^2(K-1)^2(K-2)(1+s) ^{2K-5}
	\geq 0
	\end{align*}
	and
	\begin{align*}
	\xi'''(s) \xi_0''(s) +\xi''(s) \xi_0'''(-s) =tK^2(K-1)^2(K-2)\bigl((1+s) ^{K-5}+(1+s) ^{K-2}(1-s) ^{K-3}\bigr)
	\geq 0.
	\end{align*}
	Thus, $\zeta_+$ and $\zeta_-$ are nondecreasing. The fact that $\xi_0''(s) <\xi(|s|)$ holds for all $s\in[-1,1]\setminus\{0\}$ is clear since $t\in(0,1).$ Next, assume that \eqref{eq-3} is valid. Since
	\begin{align*}
	&\xi'''(s) \xi_0''(s) -\xi''(s) \xi_0'''(s) \\
	&=K^2(K-1)^2(K-2)\bigl(t^2(1+s) ^{K-3}(1+ts) ^{K-2}-t^3(1+s) ^{K-2}(1+ts) ^{K-3}\bigr)\\
	&=(1-t)t^2K^2(K-1)^2(K-2)(1+s) ^{K-3}(1+ts) ^{K-3}\\
	&\geq 0
	\end{align*}
	and
	\begin{align*}
	&\xi'''(s) \xi_0''(-s) +\xi''(s) \xi_0'''(-s) \\
	&=K^2(K-1)^2(K-2)\bigl(t^2(1+s) ^{K-3}(1-ts) ^{K-2}+t^3(1+s) ^{K-2}(1-ts) ^{K-3}\bigr)\\
	&=(1+t)t^2K^2(K-1)^2(K-2)(1+s) ^{K-3}(1+ts) ^{K-3}\\
	&\geq 0.
	\end{align*}
	From these, $\zeta_+$ and $\zeta_-$ are nondecreasing. On the other hand, clearly
	\begin{align*}
	\xi_0''(s) =t^2K(K-1)(1+ts) ^{K-2}<K(K-1)(1+|s|)^{K-2}=\xi''(|s|)
	\end{align*}
	for all $s\in[-1,1]\setminus\{0\}.$ This completes our proof.
	\qed

\subsection{Variational representations for $\Phi_\gamma$ and $\Psi_\gamma$}

In order to establish the condition $(ii)$ in the previous section, a key ingredient we need is the variational representation for $\Phi_\gamma$ and $\Psi_\gamma$ in terms of optimal stochastic control problems. Denote by $\mathbf W=\{\mathbf{W}(w)=({W}_1(w),{W}_2(w)),\mathcal{G}_w,0\leq w\leq 1\}$ a two-dimensional Brownian motion, where the filtration $(\mathcal{G}_w)_{0\leq w\leq 1}$ satisfies the usual conditions (see Definition~2.25 in Chapter 1 of \cite{KS}). Let now $\gamma\in \mathcal{U}$ be fixed. For $0\leq s\leq 1,$ denote by $D[s]$ the space of all two-dimensional progressively measurable processes $v=(v_1,v_2)$ with respect to $(\mathcal{G}_w)_{0\leq w\leq s}$ satisfying $
\sup_{0\leq w\leq s}|v_1(w)|\leq 1$ and $\sup_{0\leq w\leq s}|v_2(w)|\leq 1.$ Endow the space $D[s]$ with the norm
\begin{align*}
\|v\|_{s}&=\Bigl(\e\int_0^s(v_1(w)^2+v_2(w)^2)dw\Bigr)^{1/2}.
\end{align*}
Recall $T(s)$ in (\ref{Tees}) and define a functional
\begin{align*}
F_\gamma^{s}(v)&=\e\left[C_\gamma^{s}(v)-L_\gamma^{s}(v)\right]
\end{align*}
for  $v\in D[s],$ where
\begin{align*}
C_\gamma^{s}(v)&:=\Psi_\gamma\Bigl(s,\int_0^s\gamma(w)T(w)v(w)dw+\int_0^sT(w)^{1/2}d\mathbf{W}(w)\Bigr),\\
L_\gamma^{s}(v)&:=\frac{1}{2}\int_0^s\gamma(w)\left<T(w)v(w),v(w)\right>dw.
\end{align*}
The functional $\Psi_\gamma$ defined above as the solution of a PDE can also be written via an optimal stochastic control problem.

\begin{proposition}\label{prop2} Let $\gamma\in\mathcal{U}.$ For any $s\in[0,1],$
	\begin{align}
	\label{prop2:eq2}
	\Psi_\gamma(0,0,0)&=\max\left\{F_\gamma^{s}(v)\,:\, v\in {D}[s]\right\}.
	\end{align}
	The maximum of \eqref{prop2:eq2} is attained by $v_{\gamma}(r)=\triangledown \Psi_\gamma(r,\vX(r))$, where the two-dimensional stochastic process $(\vX_\gamma(w))_{0\leq w\leq s}$ satisfies
	\begin{align*}
	d\vX_\gamma(w)&=\gamma(w)T(w)\triangledown \Psi_\gamma(w,\vX(w))dw+T(w)^{1/2}d\mathbf{W}(w),\\
	\vX_\gamma(0)&=(0,0).
	\end{align*}
\end{proposition}
The derivation of Proposition \ref{prop2} is a purely an application of It\^{o}'s formula. For a detailed proof, we refer the readers to \cite[Theorem 5]{C15}. Although the argument therein is for different boundary condition and $\gamma(1-)$ is bounded, the same argument carries through in the current setting with only minor modification. Note that when $q=0,$ $T$ is a diagonal matrix and 
$$
\Psi_\gamma(0,x_1,x_2)=\Phi_\gamma(0,x_1)+\Phi_\gamma(0,x_2).
$$
Proposition \ref{prop2} is a zero-temperature two-dimensional analogue of Lemma 2 in \cite{ChenAuf} (see also \cite{JagTob}). By taking $x_1=x_2=0,$ it implies the zero-temperature one-dimensional analogue of Lemma 2 in \cite{ChenAuf}, giving the following representation of $\Phi_\gamma.$ Let $W$ be a one-dimensional standard Brownian motion with respect to the filtration $(\mathcal{G}_w)_{0\leq w\leq 1}$ and $D_0[s]$ be the space of all progressively measurable processes $u$ with respect to $(\mathcal{G}_w)_{0\leq w\leq s}$ and satisfy $\sup_{0\leq w\leq s}|u(w)|dw\leq 1.$ 

\begin{corollary}\label{cor1} Let $\gamma\in \mathcal{U}$. For any $s\in[0,1],$
	\begin{align}
	\begin{split}\label{cor1:eq1}
	\Phi_\gamma(0,0)&=\max_{u\in D_0[s]} \e\Bigl[\Phi_\gamma\Bigl(s,\int_0^s\xi''(w)\gamma(w) u(w)dw+\int_0^s\sqrt{\xi''}(w)dW(w)\Bigr)\\
	&\qquad\qquad-\frac{1}{2}\int_0^s\xi''(w)\gamma(w) u(w)^2dw\Bigr].
	\end{split}
	\end{align}
    Here the maximizer is given by $u_\gamma(w)=\partial_x\Phi_\gamma(w,X_\gamma(w)),$ where $X_\gamma=(X_\gamma(w))_{0\leq w\leq s}$ is the solution to the following SDE with the initial condition $X_\gamma(0)=0$,
	\begin{align}
	\label{cor1:eq2}
	dX_\gamma&=\xi''(w)\gamma(w) \partial_x\Phi_\gamma(w,X_\gamma(w))dw+\sqrt{\xi''(w)}dW(w).
	\end{align} 

\end{corollary}

\begin{remark}\label{rmk1}
	\rm From \cite[Lemma 2]{ChenHL}, this minimizer $u_\gamma$ is unique if $\gamma(s)>0$ on $(0,1].$ In particular, from \cite[Subsection 3.2]{ChenHL}, $0$ lies in the support of $\gamma_P(ds)$, so $\gamma_P(s)>0$ on $(0,1].$
\end{remark} 
While it is generally not possible to find the solutions $\Psi_{\gamma_q}$ and $\Phi_{\gamma_P}$ explicitly and compare their values, the variational representations in Proposition \ref{prop2} and Corollary \ref{cor1} provide an elementary way to quantify the difference between $\Psi_{\gamma_q}(0,0,0)$ and $2\Phi_{\gamma_P}(0,0).$

\begin{proposition}\label{prop3}
	The following two statements hold:
	\begin{enumerate}
		\item[$(i)$] If $v_{\gamma_q}=(v_1,v_2)$ is the maximizer to the variational problem \eqref{prop2:eq2} for $\Psi_{\gamma_q}(0,0)$ using $s=|q|,$ then
		\begin{align}
		\begin{split}
		\label{sec4.3:prop:eq1}
		\Psi_{\gamma_q}(0,0,0)&\leq 2\Phi_{\gamma_P }(0,0)\\
		&-\frac{1}{2}\int_{0}^{|q|}\frac{\xi''(w)\xi_0''(\iota w)\bigl(\xi''(w)-\xi_0''(\iota w)\bigr)}{2\bigl(\xi''(w)+\xi_0''(\iota w)\bigr)^2}\gamma_P(w)\e\left(v_1(w)-\iota v_2(w)\right)^2dw.
		\end{split}
		\end{align}
		
		\item[$(ii)$] Let us define
		\begin{align}
		\begin{split}\label{sec4.3:prop:eq2}
		\left(u_1(w),u_2(w)\right)&=\frac{1}{\xi''(w)+\xi_0''(\iota w)}T(w)v_{\gamma_q}(w),
		\end{split}\\
		\begin{split}
		\label{sec4.3:prop:eq3}
		(B_1(w),{B}_2(w))&=\frac{1}{\xi''(w)^{1/2}}T(w)^{1/2}\mathbf{W}(w)
		\end{split}
		\end{align}
		for $0\leq w\leq |q|.$ If the following equality holds, $$
		\Psi_{\gamma_u}(0,0,0)=2\Phi_{\gamma_{P} }(0,0),$$
		then $u_1$ and $u_2$ are the maximizers of (\ref{cor1:eq1}) with respect to the Brownian motions $B_1$ and $B_2$ in (\ref{sec4.3:prop:eq3}) respectively, that is, on the interval $[0 ,|q|]$,
		\begin{align}
		\begin{split}
		\label{eq9}
		u_1(w)&=\partial_x\Phi_{\gamma_P }(w,X_{1,\gamma_P}(w)),\\
		u_2(w)&=\partial_x\Phi_{\gamma_P }(w,X_{2,\gamma_P}(w)),
		\end{split}
		\end{align}
		where $(X_{1,\gamma_P}(w))_{0\leq w\leq|q|}$ and $(X_{2,\gamma_P}(w))_{0\leq w\leq |q|}$ satisfy
		\begin{align}
		\begin{split}\label{eq10}
		dX_{1,\gamma_P}(w)&=\gamma_{P} (w)\xi''(w)\partial_x\Phi_{\gamma_{P} }(w,X_{1,\gamma_P}(w))dw+\xi''(w)^{1/2}dB_1(w),\\
		dX_{2,\gamma_P}(w)&=\gamma_{P}(w)\xi''(w)\partial_x\Phi_{\gamma_{P} }(w,X_{2,\gamma_P}(w))dw+\xi''(w)^{1/2}dB_2(w)
		\end{split}
		\end{align}
		with the initial condition $X_{1,\gamma_P}(0)=X_{2,\gamma_P}(0)=0.$
	\end{enumerate}
\end{proposition}  
Proposition \ref{prop3} is essentially taken from \cite[Proposition 5]{C15}. Its proof is based on the comparison between \eqref{prop2:eq2} and \eqref{cor1:eq1} with $s=|q|$. More precisely, note that 
$$
\Psi_{\gamma_q}(|q|,x_1,x_2)=\Phi_{\gamma_P}(|q|,x_1)+\Phi_{\gamma_P}(|q|,x_2).
$$ 
If one takes the optimizer $v_{\gamma_q}$ in \eqref{prop2:eq2}, then \eqref{sec4.3:prop:eq1} can be obtained after some algebraic manipulation, while the statement \eqref{eq9} follows by Remark \ref{rmk1}. For details, we refer the readers to \cite[Proposition 5]{C15}. Although the PDEs considered there have different boundary conditions, this does not affect the proof in any essential way.

\subsection{Proof of Theorem \ref{Thm2} assuming the case $(ii)$}

     First note that $\mathcal{T}_q(\gamma_q)\leq 2\mathcal{P}(\gamma_P)$ from \eqref{eq7} and \eqref{sec4.3:prop:eq1}. We claim that 
	\begin{align}
	\label{eq8}
	\mathcal{T}_q(\gamma_q)<2\mathcal{P}(\gamma_P)
	\end{align}
	for all $q\in[-1,1]\setminus\{0\}$. Assume on the contrary $\mathcal{T}_q(\gamma_q)=2\mathcal{P}(\gamma_P)$ for some $q\in[-1,1]\setminus\{0\}.$ This and \eqref{eq7} imply  $\Psi_{\gamma_q}(0,0,0)=2\Phi_{\gamma_P}(0,0)$. Note that $\gamma_P(s)>0$ on $(0,1]$ by Remark \ref{rmk1}. Since $\xi_0''(w)<\xi''(|w|)$ for all $w\in[-1,1]\setminus\{0\}$ by Lemma~\ref{lem0}, applying \eqref{sec4.3:prop:eq1} gives that $v_1=\iota v_2$ on $[0,|q|]$ and thus, the definition of $T$ in (\ref{Tees}) and \eqref{sec4.3:prop:eq2} imply that $(u_1,u_2)=\bigl(v_1,\iota v_1\bigr).$ By \eqref{eq9},
	\begin{align*}
	\partial_x\Phi_{\gamma_P}(w,X_{1,\gamma_P}(w))=u_1(w)=\iota u_2(w)=\iota\partial_x\Phi_{\gamma_P}(w,X_{2,\gamma_P}(w))
	\end{align*}
	on $[0,|q|]$. Since $\partial_x\Phi_\gamma(w,\cdot)$ is a strictly increasing odd function on $\mathbb{R}$ (see \cite[Lemma 4]{ChenHL}), $X_{1,\gamma_P}=\iota X_{2,\gamma_P}$ on $[0,|q|]$. Consequently, from \eqref{eq10},
	\begin{align*}
	0=X_{1,\gamma_P}(s)-\iota X_{2,\gamma_P}(s)&=\int_0^s\xi''(w)^{1/2}d(B_1(w)-\iota B_2(w))
	\end{align*}
	and thus, $B_1(w)=\iota B_2(w)$ on $[0,|q|]$. The definition \eqref{sec4.3:prop:eq3} then implies that 
	$$
	0=\e\bigl(B_1(w)-\iota B_2(w)\bigr)^2 =2\bigl(\xi''(w)- \xi_0''(\iota w)\bigr).
	$$
On the other hand, by Lemma \ref{lem0}, the right hand side is strictly positive -- a contradiction. 
	
	Next, from \eqref{prop2:eq2} with $s=|q|$, it is easy to see that $\Psi_{\gamma_q}(0,0,0)$ is continuous in $q.$ Therefore, by (\ref{eq8}), for any $\eps\in(0,1)$ there exists $\eta>0$ such that $\mathcal{T}_q(\gamma_q)\leq 2\mathcal{P}(\gamma_P)-3\eta$ for all $q\in S_N\setminus [-\varepsilon,\varepsilon].$ Applying (\ref{Parisi}) and \eqref{GT} yields that, for large enough $N,$ 
	\begin{align*}
	\e \MCE_N(u)\leq 2 \e \ME_N-2\eta
	\end{align*}
	for any $q\in S_N\setminus [-\varepsilon,\varepsilon].$ Furthermore, since $S_N\setminus [-\varepsilon,\varepsilon]$ contains no more than $2N$ distinct elements, using the usual concentration of measure for Gaussian extrema processes $\MCE_N(u)$ and $\ME_N$ implies that, for large enough $N,$
	\begin{align*}
	\e\MCE_N\bigl([-1,1]\setminus[-\varepsilon,\varepsilon]\bigr)\leq 2\e\ME_N-\eta.
	\end{align*}
This finishes our proof.
	\qed

\end{document}